\title{\vspace{-1cm} Independent transversals in locally sparse graphs}

\author{
Po-Shen Loh \thanks{Department of Mathematics,
Princeton University, Princeton, NJ 08544. E-mail:
ploh@math.princeton.edu.
Research supported in part by a Fannie and John Hertz Foundation Fellowship, an NSF Graduate
Research Fellowship, and a Princeton Centennial Fellowship.}
\and Benny Sudakov \thanks{Department of Mathematics,
Princeton University, Princeton, NJ 08544,  
and Institute for Advanced Study, Princeton. E-mail:
bsudakov@math.princeton.edu.
Research supported in part by NSF CAREER award DMS-0546523,
NSF
grant DMS-0355497, USA-Israeli BSF grant, Alfred P. Sloan
fellowship, and the State of New Jersey.}
}

\documentclass[11pt]{article}
\usepackage{amsmath, amssymb, epsfig, graphicx, enumerate}

\newtheorem{theorem}{Theorem}[section]

\newtheorem{lemma}[theorem]{Lemma}
\newtheorem{proposition}[theorem]{Proposition}

\newcommand{\ep}{\epsilon}
\newcommand{\pr}[1]{\mathbb{P}\left[#1\right]}
\newcommand{\E}[1]{\mathbb{E}\left[#1\right]}

\date{}
\begin{document}
\maketitle

\begin{abstract}
  Let $G$ be a graph with maximum degree $\Delta$ whose vertex set is partitioned into 
  parts $V(G)=V_1 \cup \ldots \cup V_r$. A transversal is a subset of $V(G)$ containing
  exactly one vertex from each part $V_i$.  If it is also an independent set, then we call
  it an independent transversal.  The local degree of $G$ is the maximum number of
  neighbors of a vertex $v$ in a part $V_i$, taken over all choices of $V_i$ and $v \not
  \in V_i$.  We prove that for every fixed $\epsilon > 0$, if all part sizes $|V_i| \geq
  (1+\epsilon)\Delta$ and the local degree of $G$ is $o(\Delta)$, then $G$ has an
  independent transversal for sufficiently large $\Delta$.  This extends several previous
  results and settles (in a stronger form) a conjecture of Aharoni and Holzman.  We then
  generalize this result to transversals that induce no cliques of size $s$.  (Note that
  independent transversals correspond to $s=2$.)  In that context, we prove that parts of
  size $|V_i| \geq (1+\ep)\frac{\Delta}{s-1}$ and local degree $o(\Delta)$ guarantee the
  existence of such a transversal, and we provide a construction that shows this is
  asymptotically tight.
\end{abstract}

\section{Introduction}

Let $G = (V, E)$ be a graph with maximum degree $\Delta$, whose vertices have been
partitioned into $r$ disjoint sets $V = V_1 \cup \ldots \cup V_r$.  An \emph{independent
  transversal}\/ of $G$ with respect to $\{V_i\}_{i=1}^r$ is an independent set in $G$
which contains exactly one vertex from each $V_i$. The problem of finding sufficient
conditions for the existence of an independent transversal dates back to 1975, when it was
raised by Bollob\'as, Erd\H{o}s, and Szemer\'edi \cite{BES}.  Since then, much work has been
done \cite{ABZ, AH, Al-problems-results, Al-linear-arboricity, Al-strong-chromatic,
  Ha-note, Ha-strong-chromatic, HS, Ji, Me, RS, ST, Yu}, and this basic concept has also
appeared in the study of other combinatorial problems, such as linear arboricity, strong
chromatic number and list coloring.  In particular, as part of his result on the linear
arboricity of graphs, Alon \cite{Al-linear-arboricity} used the Lov\'asz Local Lemma to
show that an independent transversal exists as long as all $|V_i| \geq 2e\Delta$.  Haxell
\cite{Ha-note} later improved his constant from $2e$ to $2$.  In the other direction, Jin
\cite{Ji} and Yuster \cite{Yu} constructed graphs with parts of size $|V_i| = 2\Delta - 1$
and maximum degree $\Delta$, with no independent transversals, but required that $\Delta$
was a power of 2.  Szab\'o and Tardos \cite{ST} recently produced constructions with the
same properties for all $\Delta$, so the constant 2 is tight.

However, in all of the above constructions, the graphs are disjoint unions of $2\Delta -
1$ complete bipartite subgraphs $K_{\Delta, \Delta}$ and the the partition into $\{V_i\}$
is done in such a way that the parts $\{V_i\}$ separate the sides of each $K_{\Delta,
  \Delta}$.  This creates many pairs of disjoint parts $(V_i, V_j)$ which have complete
bipartite subgraphs of linear size going between them.  Note that the number of edges
between such a pair $(V_i, V_j)$ is quadratic in $\Delta$.  In this paper we show that the
constant $2$ can be significantly improved if one prohibits such phenomena.  One way to
accomplish this is to introduce a constraint on the \emph{local degree}\/, which is the
maximum number of neighbors of a vertex $v$ in a part $V_i$, where $V_i$ ranges over all
parts and $v$ ranges over all vertices $v \not \in V_i$.  This constraint arises naturally
in several contexts, one of which is vertex list coloring.

Given a graph $H = (V, E)$ and a set of lists $\{C_v\}$ of available colors, one for each
vertex $v \in V$, it is a natural question to determine when we can properly color $H$
from these lists. Suppose that in addition we know that every color $c$ appears in the
lists of at most $\Delta$ neighbors of each vertex $v$; then, what minimum size lists will
guarantee a proper coloring? This question, which was proposed by Reed \cite{Re}, can be
recast as an independent transversal problem as follows.  Consider a $|V|$-partite graph
$G$ such that for each $v \in V$, $G$ has a part with $|C_v|$ vertices labeled by ordered
pairs $\{(v, c) : c \in C_v\}$.  Let two vertices $(v, c)$ and $(w,c)$ be adjacent
whenever $v$ is adjacent to $w$ in $H$ and $c \in C_v \cap C_w$.  Then $G$ has maximum
degree $\leq \Delta$ and local degree $\leq 1$, and an independent transversal in $G$ corresponds to
a proper list coloring of $H$.  (Note that not every $G$ with local degree 1 has a
corresponding list coloring problem, so this association is not reversible.)  Haxell's
result immediately implies that if all $|C_v| \geq 2\Delta$, a proper list coloring
exists.  However, this is not tight, since the local degree condition prohibits the
constructions we mentioned earlier. Indeed, for the list coloring problem Reed and Sudakov
\cite{RS} showed that in fact lists of size $(1+o(1))\Delta$ will suffice.

Aharoni and Holzman \cite{AH} adapted arguments from \cite{RS} to prove the existence of an 
independent transversal in multipartite graphs with maximum degree $\Delta$,
parts of size $(1+o(1))\Delta$, and the property that any two distinct vertices in the same 
part are at distance greater than 4 from each other.  Their result has the following 
nice application. For any collection of $n \geq (1+o(1))\Delta$
graphs $\{H_i\}_{i=1}^n$ with maximum degree $\Delta$, all
sharing the same vertex set $V$, there exists a partition $V = \bigcup_{i=1}^n I_i$ 
such that for each $i$, $I_i$ is an independent set in $H_i$.  
To see this, create a multipartite graph $G$ by making $n$ copies of each
vertex, and connect the $i$-th copy of vertex $v$ to the $i$-th copy of vertex $w$ if
$v$ is adjacent to $w$ in $H_i$. Then in $G$ there are no paths at all between any 
pair of distinct vertices in a given part. Thus we can find an independent transversal of $G$,
which gives the required partition.

Aharoni and Holzman \cite{AH} conjectured that their condition on distances could be
replaced by the weaker condition that the local degree is 1.  In this paper, we prove the
following stronger theorem, which implies their conjecture.  Our proof combines arguments
from \cite{Al-strong-chromatic} and \cite{RS}, together with some additional ideas.

\begin{theorem}
\label{main}
For every $\epsilon>0$ there exists $\gamma>0$ such that the following holds.  If $G$ is a
graph with maximum degree at most $\Delta$ whose vertex set is partitioned into parts
$V(G)=V_1 \cup \ldots \cup V_r$ of size $|V_i| \geq (1+\ep)\Delta$, and the local degree
of $G$ is at most $\gamma \Delta$, then $G$ has an independent transversal.
\end{theorem}
Note that the constant of 1 is optimal because of the following example: a disjoint union of $\Delta$
cliques of order $\Delta+1$, where each clique has exactly one vertex per part.

An independent transversal is a set with one vertex from each part $V_i$ that induces no
cliques of size 2.  Therefore, a natural generalization of this concept is the $K_s$-free
transversal, which is a transversal inducing no cliques of size $s$.  Such transversals were
recently studied by Szab\'o and Tardos \cite{ST}, who posed the problem of finding
$p(\Delta, K_s)$, which is defined to be the smallest integer $n$ that guarantees the
existence of a $K_s$-free transversal in any graph with maximum degree $\Delta$ and part
sizes at least $n$.  They provided a construction that bounds $p(\Delta, K_s) \geq
\frac{s}{s-1} \frac{\Delta}{s-1}$, and conjectured that their construction was optimal.

However, this construction also contains complete bipartite subgraphs of linear size, with
sides separated by the partition.  In light of our previous result, one may ask what can
be said when we impose a local degree restriction.  We find that we can solve that problem
asymptotically, and prove the following generalization of Theorem \ref{main}:

\begin{theorem}
  \label{generalize}
  For every $\epsilon>0$ and integer $s \geq 2$, there exists $\gamma>0$ such that the
  following holds.  If $G$ is a graph with maximum degree at most $\Delta$ whose vertex
  set is partitioned into parts $V(G)=V_1 \cup \ldots \cup V_r$ of size $|V_i| \geq
  (1+\ep)\frac{\Delta}{s-1}$, and the local degree of $G$ is at most $\gamma \Delta$, then
  $G$ has a $K_s$-free transversal.
\end{theorem}

This is asymptotically tight via a simple construction that we will give later.
Furthermore, a slight adaptation of our method proves that even without the local degree
condition, $p(\Delta, K_s) \leq 2 \big\lfloor \frac{\Delta}{s-1} \big\rfloor$, which
differs from Szab\'o and Tardos's conjecture by a factor of at most 2.  For $s=3$, this
matches their best known upper bound for $p(\Delta, K_3)$, and even is better by 1 when
$\Delta$ is odd.

The rest of this paper is organized as follows.  The next section reviews some basic
probabilistic tools we use in our proof.  In Section \ref{sec:reduce-loc-degree} we show
how to reduce Theorem \ref{main} to the special case when local degrees are bounded by a
constant. We solve this case in Section \ref{sec:const-loc-degree}.  In Section
\ref{sec:generalize}, we prove the generalization of our main result to $K_s$-free
transversals.  The final section contains some concluding remarks and open problems.
Throughout this paper we will assume wherever needed that $\gamma$ is sufficiently small.
Since, by definition, every non-trivial $r$-partite graph has local degree at least one,
this implies that $\Delta \geq \gamma^{-1}$ is sufficiently large.  We will also
systematically omit floor and ceiling signs for the sake of clarity of presentation.

\section{Probabilistic tools}
\label{sec:probabilistic-tools}

In this section we describe some classical results which we will use in our proof.  We
begin with several large-deviation inequalities.

\begin{theorem}
\label{thm:chernoff}
\label{thm:hoeffding}
(Hoeffding \cite{HMRR}, Chernoff \cite{AS}) Let $X = \sum_{i=1}^n X_i$ be a sum of bounded
independent random variables $a_i \leq X_i \leq b_i$.  Then if we let $\mu = \E{X}$,
  \begin{displaymath}
    \pr{\left|X - \mu\right| \geq t} \leq 2 \exp\left\{  -\frac{2t^2}{\sum_{i=1}^n (b_i -
      a_i)^2} \right\}.
  \end{displaymath}
In particular, when  $X_i$ are indicator variables we have
  \begin{displaymath}
    \mathbb{P}\big[| X - \mu | > t\big] < 2 e^{-2t^2/n}\,.
  \end{displaymath}
Also, for any $\ep > 0$, there exists $c_\ep >
  0$ such that
  \begin{displaymath}
    \mathbb{P}\big[| X - \mu |  > \ep \mu\big] < 2 e^{-c_\ep \mu}\,.
  \end{displaymath}
\end{theorem}

To state the next concentration result, we need to introduce two concepts.  Let $\Omega =
\prod_{i=1}^n \Omega_i$ be a probability space, and let $X : \Omega \rightarrow
\mathbb{R}$ be a random variable.
\begin{itemize}
\item Suppose that there is a constant $C$ such that changing $\omega$ in any single
  coordinate affects the value of $X(\omega)$ by at most $C$.  Then we say that $X$ is
  \emph{$C$-Lipschitz}.
\item Suppose that for every $s$ and $\omega$ such that $X(\omega) \geq s$, there exists a set $I
  \subset \{1, \ldots, n\}$ of size $|I| \leq rs$ such that every $\omega'$ that agrees
  with $\omega$ on the coordinates indexed by $I$ also has $X(\omega') \geq s$.  Then we
  say that $X$ is $r$-certifiable.
\end{itemize}

\begin{theorem}
  \label{thm:talagrand}
  (Talagrand \cite{MR}) Suppose that $X$ is a $C$-Lipschitz and $r$-certifiable random
  variable on $\Omega = \prod_{i=1}^n \Omega_i$ as above.  Then,
$$    \pr{\big|X - \E{X}\big| > t + 60C\sqrt{r\E{X}}} \leq 4e^{ 
-\frac{t^2}{8C^2r\E{X}}}\,.$$
\end{theorem}

Finally we need the symmetric version of the Lov\'asz Local Lemma, which is typically used
to show that with positive probability, no ``bad'' events happen.
\begin{theorem}
  (Lov\'asz Local Lemma \cite{AS}) Let $A_1, \ldots, A_n$ be events 
in a probability space. Suppose that there exist constants $p$ and $d$ such that all 
$\pr{A_i} \leq  p$, and  each event $A_j$ is mutually independent of all of the other 
events $\{A_i\}$ except at most $d$ of them.  If $ep(d+1) \leq 1$, where $e$ is the base of 
the natural logarithm, then  $\pr{\bigcap \overline{A_i}} > 0$.
\end{theorem}

The following result can be deduced quickly from this lemma. We record it here for later
use, and sketch the proof for the sake of completeness.
\begin{proposition}
  \label{prop:2e} (Alon \cite{Al-linear-arboricity})
  Let $G$ be a multipartite graph with maximum degree $\Delta$, whose parts $V_1, \ldots,
  V_r$ all have size $|V_i| \geq 2e\Delta$.  Then $G$ has an independent transversal.
\end{proposition}
The proof follows by applying the Local Lemma to the probability space where we
independently and uniformly select one vertex from each $V_i$.  For each edge $e$ of $G$,
let the ``bad'' event $A_e$ be when both endpoints of $e$ are selected.  The dependency is
bounded by $2 (2e\Delta) \Delta -1$, and the probability of each $A_e$ is at most
$(2e\Delta)^{-2}$, so the Local Lemma implies this statement immediately.

\section{Reducing local degrees}
\label{sec:reduce-loc-degree}

In this section, we show that it is enough to prove Theorem \ref{main} only in the case
when the local degree is bounded by a constant. This will be an immediate consequence of
the following claim.

\begin{theorem}
\label{t21}
For any $\ep > 0$, there exists $\gamma_0 > 0$ such that for all $\gamma < \gamma_0$ and all $\Delta$, the following holds.
Let $G$ be a
multipartite graph with maximum degree $\leq \Delta$, parts $V_1, \ldots, V_r$ of size $|V_i|
\geq (1+\ep)\Delta$, and local degree $\leq \gamma \Delta$.  Then there exist subsets
$W_i \subset V_i, 1 \leq i \leq r$ such that the $r$-partite subgraph $G'$ of $G$ induced
by the set $\bigcup W_i$ has the following properties. The maximum degree of $G'$ is at
most some $\Delta' > \gamma^{-1/3}$, each $W_i$ has size $\geq (1+\ep/8)\Delta'$ and
the local degree of $G'$ is less than 10.
\end{theorem}

We first prove the following special case of the above theorem, when $\Delta^{2/3} \leq
\gamma^{-1}$.

\begin{lemma}
\label{l21}
For any $0 < \ep < 1$, there exists $\Delta_0$ such that the following holds for all $\Delta > \Delta_0$.
Let $G$ be a multipartite graph with maximum degree 
$\leq \Delta$, parts $V_1, \ldots, V_r$ of size $|V_i| \geq (1+\ep)\Delta$, and local degree 
$\leq \Delta^{1/3}$. Then there exist subsets $W_i
\subset V_i, 1 \leq i \leq r$ such that the $r$-partite subgraph $G'$ of $G$ induced by
the set $\bigcup W_i$ has the following properties. The maximum degree of $G'$ is 
at most $\Delta' = (1+\ep/3)\Delta^{1/3}$, each $W_i$ has size at least $(1+\ep/4)\Delta'$ and
the local degree of $G'$ is less than 10.
\end{lemma}

\noindent {\bf Proof.}\, By discarding vertices, we may assume that all $|V_i| = (1+\ep)\Delta$.  For every $1 \leq i \leq r$, create $W_i$ by
choosing each vertex of $V_i$ randomly and independently with probability $p=\Delta^{-2/3}$.  
Define the following three types of bad events.  For each vertex $v$, let $A_v$ be the event 
that the number of neighbors of $v$ in $W=\bigcup W_j$ exceeds $(1+\ep/3)\Delta^{1/3}$.  
For each vertex $v$ and part $V_i$ in which $v$ has at least one neighbor, let $B_{v,i}$ be
the event that the number of neighbors of $v$ in $W_i$ is at least $10$.  Finally, for every 
$1
\leq i \leq r$, let $C_i$ be the event that $|W_i| < (1+2\ep/3)\Delta^{1/3}$.  Note that
$(1+2\ep/3)\Delta^{1/3} = \frac{1+2\ep/3}{1+\ep/3}\Delta'$, which exceeds
$(1+\ep/4)\Delta'$ if $\ep < 1$ (which we assumed).  We use the symmetric version of
the Local Lemma to show that with positive probability, no bad events happen.

To bound the dependency, observe that each of the events $A_v$, $B_{v,i}$ is completely determined by the choices
on all vertices within distance one from $v$, and $C_k$ is completely determined by the choices on all vertices in $V_k$.
Since degrees are bounded by $O(\Delta)$ and all $|V_k| \leq O(\Delta)$, each event is mutually independent of all but
$O(\Delta^2)$ events.

Now we compute the probabilities of bad events. Since the number of neighbors of a vertex
$v$ in $W$ is binomially distributed with mean at most $\Delta p =\Delta^{1/3}$, the
standard Chernoff estimate (Theorem \ref{thm:chernoff}) implies that the probability that
it exceeds $(1+\ep/3)\Delta^{1/3}$ is at most $e^{-\Omega(\Delta^{1/3})}\ll \Delta^{-3}$.
Similarly, the size of the set $W_i$ is binomially distributed with mean at least
$(1+\ep)\Delta^{1/3}$. Hence, using the Chernoff estimate again, we conclude that $\pr{C_i}
\leq e^{-\Omega(\Delta^{1/3})}\ll \Delta^{-3}$. Finally, since the number of neighbors of
$v$ in $V_i$ is bounded by $\Delta^{1/3}$, we have
\begin{displaymath}
  \pr{B_{v,i}}\leq {\Delta^{1/3} \choose 10}p^{10} \leq \Delta^{-10/3} \ll \Delta^{-3}
\end{displaymath}
Thus, by the Local Lemma,  with positive probability none of the events
$A_v$, $B_{v,i}$ and $C_i$ happen and we obtain an induced subgraph $G'$ of $G$ 
which has all the desired properties. \hfill $\Box$  

The general case of Theorem \ref{t21} cannot be proved using the above arguments, since if
$\gamma^{-1}$ were much smaller than $\log \Delta$, the number of dependencies would
overwhelm the probabilities in the application of the Local Lemma.  To overcome this
difficulty, we follow an approach similar to the one used in \cite{Al-strong-chromatic}
and construct the desired subgraph by a sequence of random halving steps. This is done
via the following lemma.

\begin{lemma}
\label{alon}
Let $G$ be a multipartite graph with maximum degree at most $\Delta$, parts $V_1, \ldots,
V_r$ each of size $2s$, and local degree at most $d$. Suppose that $\Delta$ is
sufficiently large and $d>\log^4 \Delta$. Then there exist subsets $U_i \subset V_i$, each
of size $s$, such that the subgraph of $G$ induced by $\bigcup U_i$ has maximum degree at
most $\Delta/2 + \Delta^{2/3}$, and local degree at most $d/2 + d^{2/3}$.
\end{lemma}

\noindent {\bf Proof.}\, Within each $V_i$, arbitrarily pair up the vertices so that each
vertex $v$ has a mate $M(v)$. Note that this pairing doesn't need to have any correlation
with the original edges of $G$.  For each pair of vertices $\{v, M(v)\}$, randomly and
independently designate one of the vertices to be in $U_i$.  Clearly all $U_i$ will have
size $s$. For each vertex $v$, let $A_v$ be the event that the number of neighbors of $v$
in $U = \bigcup U_i$ exceeds $\Delta/2 + \Delta^{2/3}$. Also for every part $V_k$ and
vertex $v \not \in V_k$, let $B_{v, k}$ be the event that the number of neighbors of $v$
in $U_k$ exceeds $d/2 + d^{2/3}$.  We will use the Local Lemma again to prove that with
positive probability none of these events occurs.

Fix a vertex $v$ and consider the event $A_v$. Note that if two neighbors of $v$ are
paired with each other by $M$, then exactly one of them will be in $U$. Let $T$ be the set
of all neighbors of $v$ which are paired by $M$ to vertices which are not neighbors of
$v$. Then the number of neighbors of $v$ in $U$ is at most $(\Delta-|T|)/2$ plus the
number of members of $T$ that belong to $U$.  The second number is binomially distributed
with parameters $|T| \leq \Delta$ and $1/2$. Therefore by the Chernoff bound (Theorem
\ref{thm:chernoff}), we have that the probability that it deviates from its mean by at
least $\Delta^{2/3}$ is bounded by $2e^{-2(\Delta^{2/3})^2/|T|}\ll \Delta^{-3}$. Using
similar arguments, together with the assumption that $d>\log^4 \Delta$, we can bound
$\pr{B_{v, k}} \leq 2e^{-2(d^{2/3})^2/d} \ll \Delta^{-3}$.

To bound the dependency, observe that we can argue exactly as in the proof of the previous
lemma to show that every bad event depends on at most $O\big(\Delta^2\big)$ other such
events. Thus, by the Local Lemma we have that with positive probability none of the events
$A_v$, $B_{v, k}$ happen.  \hfill $\Box$

\vspace{0.15cm}
\noindent {\bf Proof of Theorem \ref{t21}.}\, Let $G$ be a
multipartite graph with maximum degree $\leq \Delta$, local degree $\leq
d=\gamma\Delta$, and parts $V_1, \ldots, V_r$ of size
$(1+\ep)\Delta$. First consider the case when $\gamma^{-4/3}\geq
\Delta$. Then $d \leq \Delta^{1/4}$, and the result of the theorem
follows from Lemma \ref{l21} because $\Delta' > \Delta^{1/3} \geq
\gamma^{-1/3}$, since $\Delta \geq \gamma^{-1}$ as was noted at the end
of the introduction.

It remains to consider the case $\gamma^{-4/3}< \Delta$.  (We choose $-4/3$ because then our argument will give $\Delta' > \gamma^{-1/3}$.)  Let $j\geq 1$ be the integer for
which $2^{j-1}<\gamma^{4/3}\Delta\leq 2^j$. By deleting at most $2^j$ vertices from each
$V_i$, we may assume that the size $n > (1+\ep)\Delta-2^j$ of every part is divisible by
$2^j$.  Define the sequences $\{\Delta_t\}$ and $\{d_t\}$ by setting $\Delta_0=\Delta$,
$d_0=d=\gamma\Delta$, and 
\begin{displaymath}
  \Delta_{t+1}=\frac{\Delta_t}{2}+\Delta_t^{2/3}, \hspace{1.5cm} d_{t+1}=\frac{d_t}{2}+d_t^{2/3}\, .
\end{displaymath}
We claim that: 
$$ {\bf (i)}\,\, \gamma^{-4/3}/2 < \Delta_j \leq (1+\ep/4)\frac{\Delta}{2^j}, \quad \quad 
\quad  {\bf (ii)}\,\, d_j \leq 8\Delta_j^{1/4},  \quad \quad \quad
{\bf (iii)}\,\, d_t > \log^4 \Delta_t ~~\forall~ 0 
\leq t < 
j.$$

Suppose this is true.  By (iii), we can apply Lemma \ref{alon} to split each part
$V_i$ in half and obtain a new $r$-partite graph $G_1$ with maximum degree at most
$\Delta_1$ and local degree at most $d_1$.  Continuing in this manner for $j$ iterations,
applying Lemma \ref{alon} to split the graph in half each time, we obtain a sequence of
$r$-partite graphs $G \supset G_1 \cdots \supset G_j$.  Note that $\Delta_t$ and $d_t$
are upper bounds on the maximum and local degrees of each $G_t$, respectively.  Moreover,
all parts in each $G_t$ have size $n_t=n/2^t$.

By the lower bound of (i), we can make $\Delta_j$ as large as necessary by decreasing
$\gamma$, so the upper bound of (i) yields $n_j > \frac{(1+\ep)\Delta-2^j}{2^j} \geq
\frac{1+\ep}{1+\ep/4}\Delta_j - 1 > (1+\ep/2)\Delta_j$ (assume $\ep < 1$).  By
(ii), $d_j \ll \Delta_j^{1/3}$.  Applying Lemma \ref{l21} to $G_j$ with
$\ep/2$ instead of $\ep$, we obtain a new subgraph $G'$ with maximum degree at most
$\Delta' = (1+\ep/6)\Delta_j^{1/3} > (\gamma^{-4/3}/2)^{1/3} \gg
\gamma^{-1/3}$, part sizes at least $(1+\ep/8)\Delta'$, and local degree less than 10.
This completes the proof of the theorem.

To finish we need to prove our claim.  The lower bound of (i) follows immediately
from the definition of $j$, because $\Delta_j \geq \Delta/2^j > \gamma^{-4/3}/2$.  For the
upper bound, $\Delta_{t+1}=\Delta_t/2+\Delta_t^{2/3}\leq
\frac{1}{2}\big(\Delta_t^{1/3}+1)^3$, so taking cubic roots and subtracting
$1/(2^{1/3}-1)$ from both sides, we obtain
\begin{displaymath}
  \Delta_{t+1}^{1/3}-\frac{1}{2^{1/3}-1}
\leq \frac{1}{2^{1/3}}\big(\Delta_t^{1/3}+1\big)- \frac{1}{2^{1/3}-1}=
\frac{1}{2^{1/3}}\left(\Delta_t^{1/3}-\frac{1}{2^{1/3}-1}\right).
\end{displaymath}
Therefore,
\begin{displaymath}
  \Delta_j^{1/3}-\frac{1}{2^{1/3}-1}
\leq \frac{1}{2^{j/3}} \left(\Delta_0^{1/3}-\frac{1}{2^{1/3}-1}\right),
\end{displaymath}
and since $\Delta_0=\Delta$ and $2^{1/3}-1>1/4$,
\begin{displaymath}
  \Delta_j^{1/3} \leq \frac{\Delta^{1/3}}{2^{j/3}}+4
\leq (1+\ep/4)^{1/3}\frac{\Delta^{1/3}}{2^{j/3}}.
\end{displaymath}
The last inequality follows from our assumption that $\gamma$ is small and hence
$\Delta/2^j>\gamma^{-4/3}/2$ is large enough.  Therefore $\Delta_j \leq
(1+\ep/4)\Delta/2^j$. Note that since $\Delta_j \geq \Delta_t/2^{j-t}$, we have $\Delta_t
\leq (1+\ep/4) \Delta/2^t$ for all $t<j$ (we will use this in the proof of (iii)).

For (ii), the same argument as above (just substitute $d_t$ for $\Delta_t$) shows
that
\begin{displaymath}
  d_j^{1/3} \leq \frac{d^{1/3}}{2^{j/3}}+4\leq 2\frac{d^{1/3}}{2^{j/3}},
\end{displaymath}
where the last inequality used that $d/2^j = \gamma \Delta/2^j > \gamma^{-1/3}/2$ is
large.  Hence $d_j\leq 8d/2^j = 8 \gamma \Delta/2^j$.  By definition of $j$ and $\{\Delta_t\}$,
$\gamma^{4/3} \Delta \leq 2^j$, so $d_j \leq 8 \big(\Delta/2^j\big)^{1/4}\leq
8\Delta_j^{1/4}$.

Finally, to prove (iii), note that by definition of $j$, 
$\gamma >(\Delta/2^{t})^{-3/4}$ for all $t<j$. Thus
\begin{displaymath}
  d_t \geq d/2^t =\gamma \Delta/2^{t} \geq
  (\Delta/2^{t})^{-3/4} (\Delta/2^{t})=(\Delta/2^{t})^{1/4} \geq
  \left(\frac{\Delta_t}{1+\ep/4}\right)^{1/4} \gg \log^4 \Delta_t,
\end{displaymath}
and we are done. \hfill $\Box$

\section{Transversals in graphs with constant local degree}
\label{sec:const-loc-degree}
In this section we obtain the following result, which completes the 
proof of our main theorem.

\begin{theorem}
  For any $\ep > 0$ and constant $C$, the following holds for all sufficiently large
  $\Delta$.  Let $G$ be a multipartite graph with maximum degree $\leq \Delta$, parts
  $V_1, \ldots, V_r$ of size $|V_i| \geq (1+\ep)\Delta$, and local degree $\leq C$.
  Then $G$ has an independent transversal.
\end{theorem}

The proof of this result is based on the approach from \cite{RS} together with some additional 
ideas. We use the semi-random method, which constructs an independent transversal in several 
iterations.  Each iteration is a random procedure, for which we prove that there is a choice of random bits which give desirable output.  
We then fix that choice and assume it as the state of affairs for the next iteration.  Consider the following random process, 
which will provide us with one iteration of our algorithm.
\begin{enumerate}
\item Activate (for this iteration) each remaining part independently
  with probability $1/\log \Delta$.
\item Uniformly at random select a vertex from each activated part and denote by
$T$ the set of all selected vertices.
\item For each $i$ and $v \in V_i \cap T$, if $v$ is not adjacent to any $w \in V_j \cap
  T$ with $j < i$, then add $v$ to the independent transversal.
\item For each vertex $v$ added to the independent transversal in Step 3, delete the entire part containing it from $G$.
  Also delete all neighbors of all vertices in $T$ from $G$.
\end{enumerate}
Observe that the deletions ensure that after each iteration, the partial
independent transversal constructed so far has no adjacencies among the remaining vertices.  Our objective
will be to show that after performing several iterations, the remaining graph will have
maximum degree $\leq \Delta'$ and parts of size $\geq 2e\Delta'$, for some $\Delta'$.
Then, we will abort the algorithm, and apply Proposition \ref{prop:2e} to complete the
construction of our independent transversal in a single step.

\subsection{Setting the stage}

In our study of the evolution of degrees and part sizes, the following definitions are
useful.  For each part $V_i$, let $s_t(i)$ be its size at the start of iteration $t$.  For
each vertex $v$, let $N_t(v)$ be the set of $v$'s neighbors at the start of iteration $t$,
and let $d_t(v) = |N_t(v)|$.

Next, define the sequences $\{S_t\}$ and $\{D_t\}$ by setting $S_1 = (1+\ep)\Delta$, 
$D_1 = \Delta$, and 
\begin{displaymath}
  S_{t+1} = \left( 1 - \frac{1}{(1 + 3\ep/4) \log \Delta} \right) S_t\, , \hspace{1.5cm}
  D_{t+1} = \left( 1 - \frac{1}{(1 + \ep/4) \log \Delta} \right) D_t\, .
\end{displaymath}
Let $\mathbf{P}(t)$ be the property that at the start of iteration $t$, all remaining
parts have size at least $S_t$, and all remaining vertices $v$ have $d_t(v) \leq D_t$.
(Completely ignore deleted parts and vertices.)  We will prove by induction that there is always a choice of random bits such that  we can perform
iterations with property
$\mathbf{P}(t)$ holding for every $t
\leq 2 + \frac{10}{\ep} \log \Delta$. Then at the end of iteration $t' = \lceil
\frac{10}{\ep} \log \Delta \rceil$, all remaining parts have size at least $S_{t'+1}$ and
all remaining vertices have degree at most $D_{t'+1}$.  A routine calculation reveals that
\begin{eqnarray*}
  \frac{D_{t' + 1}}{S_{t'+1}} &=& \frac{\left(1 - \frac{1}{(1 + \ep/4)\log
        \Delta}\right)^{t'} D_1}{\left(1 - \frac{1}{(1 + 3\ep/4)\log \Delta}\right)^{t'}
    S_1} \ \ \leq \ \ \left(\frac{1 - \frac{1}{(1 + \ep/4)\log \Delta}}{1 - \frac{1}{(1 + 3\ep/4)\log
        \Delta}} \right)^{t'} \\
  &\leq& \left(1 - \frac{\ep}{5 \log \Delta}\right)^{\frac{10 \log \Delta}{\ep}} \ \ \leq
  \ \  e^{-2} \ \ < \ \  \frac{1}{2e}.
\end{eqnarray*}
Therefore, by Proposition \ref{prop:2e} there is an independent transversal through the
remaining parts, as promised above.  This will have no adjacencies with the partial independent
transversal constructed by the first $t'$ iterations, so the union of the two partial
transversals will be a full independent transversal.  Note that if $t \leq 1 +
\frac{10}{\ep} \log \Delta$, then $S_t = \Theta(\Delta)$ and $D_t = \Theta(\Delta)$. We
will use this fact throughout the rest of the proof.

It remains to show that if at the beginning of iteration $t$ we have a graph with property
$\mathbf{P}(t)$, then with positive probability the graph obtained at the end of this
round satisfies $\mathbf{P}(t+1)$. Define the following family of bad events.  Let $A_i$
be the event that $s_{t+1}(i) < S_{t+1}$ and let $B_v$ be the event that $d_{t+1}(v) >
D_{t+1}$.  The dependencies among these events are polynomial in $\Delta$. To see this
consider the auxiliary graph $H$ obtained by adding edges such that every part $V_i$
becomes a clique.  If we know the algorithm's choices on the ``patch'' consisting
of all vertices within distance (with respect to edges in $H$) 4 from
$v$, then $B_v$ is completely determined.  
This is because a neighbor $w$ of $v$ can only be deleted in two ways: either a neighbor of $w$ is selected in Step 2,
or the entire part containing $w$ is deleted because a vertex $x$ in that part is selected, but none of $x$'s neighbors
in lower-indexed parts are selected.
So, each event $B_v$ is mutually independent from all other events $B_w$ whose patches are disjoint from its own.
Since the part sizes are $O(\Delta)$, the degrees in $H$ are also $O(\Delta)$, so the dependency is bounded by $O(\Delta^8)$.
Events of
type $A_i$ are determined by even smaller patches, so the total dependency is 
also
$O(\Delta^8)$. Therefore if we prove that for every part $V_i$ and vertex $v$
$${\bf (i)}\,\, \pr{s_{t+1}(i) < S_{t+1}} \ll  e^{-\log \Delta 
\log \log  \Delta}     \quad\quad
\mbox{and} \quad\quad {\bf (ii)}\,\, \pr{d_{i+1}(v) > D_{i+1}} \ll e^{-\log \Delta 
\log \log
  \Delta},$$
then we can apply the Local Lemma to show that with positive probability none of 
the events $A_i$, $B_v$ occur. This corresponds precisely to property
$\mathbf{P}(t+1)$, completing the induction step.  Thus it remains to establish the
two probability bounds above.

\subsection{Parts remain large enough}
Suppose that our graph has
property $\mathbf{P}(t)$, and let $V_i$ be some part of this graph. In this section we bound the 
probability that the size of $V_i$ at the end of iteration $t$ is less than $S_{t+1}$. 

For every vertex $v$ and  part $V_k$, define
$d_t(v, k)$ to be the number of neighbors of $v$ in part $V_k$ at the start of iteration
$t$.  Since $D_t/S_t < D_1/S_1 = 1/(1+\ep)$, by linearity of expectation we have

\begin{eqnarray*}
  \E{s_{t+1}(i)} &=& \sum_{v \in V_i} \prod_{k=1}^r \left( 1 - \frac{1}{\log \Delta}
    \frac{d_t(v, k)}{s_t(k)} \right) \\
&\geq&  \sum_{v \in V_i} \left( 1 - \frac{1}{\log \Delta}
    \frac{\sum_kd_t(v,k)}{S_t} \right) \\
&=& \sum_{v \in V_i} \left( 1 - \frac{1}{\log \Delta}
    \frac{d_t(v)}{S_t} \right) \\
  &\geq& s_t(i) \left( 1 - \frac{1}{\log \Delta}
    \frac{D_t}{S_t} \right) \\
 &>& s_t(i) \left( 1 - \frac{1}{\log \Delta}
    \frac{1}{1+\ep} \right).
\end{eqnarray*}
Instead of proving concentration of $s_{t+1}(i)$ directly, we consider the number of 
vertices we deleted from the part $V_i$ in the $t$-th iteration and prove that this 
random variable $R =
s_t(i) - s_{t+1}(i)$ is concentrated.  Since the local degree is bounded by $C$, changing the assignment of 
any vertex can change $R$ by at most $C$. Also, if $R\geq s$, there are at most $s$ vertices in $T$, 
each with neighbor(s) in $V_i$, such that their selection certifies that $R\geq s$. 
Therefore $R$ is $C$-Lipschitz and 1-certifiable. Note that
$R \leq s_t(i) = \Theta(\Delta)$ and $\sqrt{\E{R}} \ll s_t(i) /\log^2 \Delta$. Thus,
using Talagrand's inequality (Theorem \ref{thm:talagrand}), we obtain 

\begin{displaymath}
  \pr{|R - \E{R}| > \frac{s_t(i)}{\log^2 \Delta}} < \exp\left\{ -\Theta\left(
      \frac{s_t(i)}{\log^4 \Delta}  \right) \right\} \ll e^{-\log \Delta \log \log \Delta}.
\end{displaymath}
Now for
sufficiently large $\Delta$,
\begin{eqnarray*}
  S_{t+1} &\leq& \left(1 - \frac{1}{(1+3\ep/4)\log \Delta} \right) s_t(i)\\ 
  &\leq& \left(1 - \frac{1}{(1+\ep)\log \Delta} - \frac{1}{\log^2
        \Delta} \right) s_t(i)\\ 
  &\leq& \E{s_{t+1}(i)} - \frac{s_t(i)}{\log^2 \Delta}.
\end{eqnarray*}
Note that since we fixed the output of the $(t-1)$-st iteration, the value of $s_t(i)$ in the definition of $R$ is fixed as well.
Thus
by linearity of expectation, $s_{t+1}(i) - \E{s_{t+1}(i)} = \E{R} - R$, so
\begin{eqnarray*}
\pr{s_{t+1}(i) < S_{t+1}} &\leq& \pr{s_{t+1}(i) <
\E{s_{t+1}(i)} - \frac{s_t(i)}{\log^2 \Delta}}\\
&\leq& \pr{\Big|s_{t+1}(i)-\E{s_{t+1}(i)}\Big| > \frac{s_t(i)}{\log^2 \Delta}}\\
&=& \pr{\big|R - \E{R}\big| > \frac{s_t(i)}{\log^2
      \Delta}} \ll e^{-\log \Delta \log \log \Delta}\, ,
\end{eqnarray*}
which implies (i).

\subsection{Degrees shrink quickly enough}
In this section we prove that if our graph has 
property $\mathbf{P}(t)$ then for every vertex $v$ the probability that 
its degree at the end of iteration $t$ is greater than $D_{t+1}$ is 
$\ll e^{-\log \Delta \log \log \Delta}$. Fix  a vertex $v$.  If we have $d_t(v) \leq D_{t+1}$, then
we are already done, so suppose that is not the case.  Then $\Delta \geq d_t(v) > D_{t+1} =
\Theta(\Delta)$.  For each vertex $v$, let $z_t(v)$ be the number of neighbors of $v$
whose entire part was deleted in Step 4 of iteration $t$.  Clearly $d_{t+1}(v)
\leq d_t(v) - z_t(v)$, so if $z_t(v) \geq \frac{d_t(v)}{\log \Delta} -
\Theta\left(\frac{d_t(v)}{\log^2 \Delta}\right)$, then for sufficiently large $\Delta$ we have
\begin{eqnarray*}
  d_{t+1}(v) \ \leq \ d_t(v) - z_t(v) &\leq& \left[1 - \frac{1}{\log \Delta} - \Theta\left(\frac{1}{\log^2
        \Delta}\right)\right] d_t(v) \\
  &\leq& \left[1 - \frac{1}{(1+\ep/4)\log \Delta} \right] D_t \ = \ D_{t+1}.
\end{eqnarray*}
Thus to prove (ii) it is enough to show
\begin{equation}
  \label{z-goal}
  \pr{z_t(v) < \frac{d_t(v)}{\log \Delta} - \Theta\left( \frac{d_t(v)}{\log^2 \Delta}
    \right)} \ll e^{-\log \Delta \log \log \Delta}.
\end{equation}

Recall our notation that for a vertex $v$ and a part $V_k$, $d_t(v, k)$ is the number of
neighbors of $v$ in $V_k$.  Call a part $V_k$ \emph{relevant} for $v$ if $d_t(v, k) \geq 1$,
i.e., $v$ has at least one neighbor in $V_k$. To analyze the behavior of $z_t(v)$, we divide 
the $t$-th iteration of the algorithm into 2 independent phases.  

\vspace{0.1cm}
\noindent
{\bf Phase I.}\, Activate each part relevant for $v$ independently with probability
  $1/\log \Delta$, and define the random variable
  \begin{displaymath}
    X_1 = \sum_{k=1}^r d_t(v, k) I_1(k)\,,
  \end{displaymath}
where the indicator $I_1(k)=1$ if part $V_k$ was activated and zero otherwise.
 Randomly select a vertex from each of these activated parts, and collect the selected
  vertices in a set $T_1$.  Define the subset $S \subseteq T_1$ as follows.  
For every $i$ and $x \in V_i \cap T_1$, we put it in $S$ iff  $x$ is not adjacent to any $y \in V_j \cap 
T_1$ with $j < i$. Let $I_2(k)$ be an indicator random variable which equals one iff 
$V_k \cap S \not =\emptyset$, and define
  \begin{displaymath}
    X_2 = \sum_{k=1}^r d_t(v, k) I_2(k)\,.
  \end{displaymath}

\noindent
{\bf Phase II.}\, Activate the rest of the parts (i.e., parts that are not relevant for
$v$) independently, each 
with  probability $1/\log \Delta$, and randomly select a vertex from each of them.  Let
$T_2$ be the set of vertices selected in this phase.  For each
  $i$ and $u \in V_i \cap S$, if $u$ is adjacent to some $w \in V_j \cap T_2$ with $j < i$,
  then remove $u$ from $S$.  Define the random variable
  \begin{displaymath}
    X_3 = \sum_{k=1}^r d_t(v, k) I_3(k)\,,
  \end{displaymath}
where the indicator $I_3(k)=1$ iff part $V_k$ still has at least one vertex in $S$.
\vspace{0.1cm}

Observe that, by definition, the parts relevant for $v$ which we delete entirely during iteration $t$ are  
exactly the ones with $I_3(k)=1$. Therefore $z_t(v) = X_3 \leq X_2 \leq X_1$.  
Our strategy will be to bound $z_t(v)$ by
starting from $X_1$ and working towards $X_3$.  By linearity of expectation, 
$\E{X_1} = \sum_k \frac{d_t(v,k)}{\log \Delta}=\frac{d_t(v)}{\log \Delta}$.  Also, since 
$d_t(v) = \Theta(\Delta)$ (see the beginning of this section), local degrees are $\leq C$, and the number of
nonzero $d_t(v, k)$ is at most $\Delta$, we can apply Hoeffding's inequality (Theorem
\ref{thm:hoeffding}) to the sum of
the terms in $X_1$ with $d_t(v, k) \neq 0$ and conclude that
\begin{equation}
  \label{x1}
  \pr{|X_1 - \E{X_1}| > \frac{d_t(v)}{\log^2 \Delta}} \leq 2 \exp\left\{
    -\frac{2}{\Delta C^2} \left(\frac{d_t(v)}{\log^2 \Delta}\right)^2  \right\} \ll
  e^{-\log \Delta \log \log \Delta}.
\end{equation}

Next, let us estimate $X_2$ by studying the difference $X_1 - X_2$.  Reveal the random
selections in the parts  activated in Phase I in order of part number (i.e. if $i<j$ and $V_i$ and
$V_j$ were activated, reveal the vertex selection in $V_i$ first).  For each activated
part $V_i$, the difference $X_1 - X_2$ will gain $d_t(v, i)$ precisely when the selected
vertex $x \in V_i \cap T_1$ is adjacent to some selected vertex $y \in V_j \cap T_1$ with
$j < i$.  Call such an event a \emph{conflict}. Its probability is at most $\frac{C
  X_1}{S_t}$, because there are at most $X_1$ activated parts with $j < i$, each of their
selected vertices has degree at most $C$ into $V_i$, and $|V_i| \geq S_t$ by 
property $\mathbf{P}(t)$.  Now condition on $|X_1 - \E{X_1}| \leq
\frac{d_t(v)}{\log^2 \Delta}$.
If $N \leq X_1$ is the number of parts activated in Phase I, the probability that
there are $\geq 4C \frac{d_t(v)}{\log^2 \Delta}$ conflicts is bounded by
\begin{eqnarray*}
  {N \choose 4C\frac{d_t(v)}{\log^2 \Delta}}
  \left(\frac{CX_1}{S_t}\right)^{4C\frac{d_t(v)}{\log^2 \Delta}} &\leq& 
  \left(\frac{e X_1}{4C \frac{d_t(v)}{\log^2 \Delta}} \frac{C X_1}{S_t}
  \right)^{4C\frac{d_t(v)}{\log^2 \Delta}} \\
  &\leq& \left(
    \frac{e}{4} \frac{ \big(\frac{d_t(v)}{\log \Delta} + \frac{d_t(v)}{\log^2
        \Delta}\big)^2 }{\frac{d_t(v)}{\log^2 \Delta} \, S_t}
  \right)^{4C\frac{d_t(v)}{\log^2 \Delta}} \\
  &\leq& \left( \frac{e+0.1}{4} \right)^{4C\frac{d_t(v)}{\log^2 \Delta}} \\
  &\ll& e^{-\log \Delta \log \log \Delta}.
\end{eqnarray*}
Here we used that $S_t \geq d_t(v)$ and $\Delta$ is sufficiently large.
Since all $d_t(v, i) \leq C$, each conflict can account
for a value gain of at most $C$ in $X_1-X_2$. Therefore, we proved that conditioned
on $|X_1 - \E{X_1}| \leq \frac{d_t(v)}{\log^2 \Delta}$,
\begin{equation}
  \label{x2}
  \pr{X_1 - X_2 \geq 4C^2 \frac{d_t(v)}{\log^2 \Delta}} \ll e^{-\log \Delta \log \log \Delta}.
\end{equation}

To estimate $X_3$, we will use Talagrand's Inequality (Theorem \ref{thm:talagrand}) to show
that the difference $X_2-X_3$ is strongly concentrated.  This requires a Lipschitz
condition, so let us first ensure that we have a good Lipschitz constant.  Let $W$ be the set of vertices
which have at least one neighbor in some part relevant for $v$. Since there are 
at most $D_t$ parts relevant for $v$, it is easy to see that $|W| \leq D_t^2S_t\leq (1+\ep)\Delta^3$. 
For $w \in W$, let $B_w$ be 
the event that at least $\log \Delta$ neighbors of $w$ are selected for $T_1$ in Phase I. The number of neighbors of $w$
in a given part is at most $C$, so the probability that one of them appears in $T_1$ is 
$\leq \frac{C}{S_t \log \Delta}$, and this happens independently for distinct parts. Since
$w$ has neighbors in at most $\Delta$ parts and 
$S_t=\Theta(\Delta)$, we obtain
$$\pr{B_w} \leq {\Delta \choose \log \Delta} \left(
\frac{C}{S_t \log \Delta}\right)^{\log \Delta} \leq 
\left( \frac{e \Delta}{\log \Delta} \frac{C}{S_t \log \Delta}\right)^{\log \Delta} 
\ll e^{-1.5\log \Delta \log \log \Delta}.$$
This implies that
\begin{equation}
\label{lipschitz}
\pr{\bigcup B_w} \leq (1+\ep)\Delta^3 e^{-1.5\log \Delta \log \log \Delta}
\ll e^{-\log \Delta \log \log \Delta}.
\end{equation}
Combining inequalities \eqref{x1}, \eqref{x2}, and
\eqref{lipschitz}, we see that
\begin{equation}
  \label{ready-for-x3}
  \pr{\left\{\frac{d_t(v)}{\log \Delta} -
      5C^2\frac{d_t(v)}{\log^2 \Delta} \leq X_2 \leq \frac{d_t(v)}{\log \Delta} +
      \frac{d_t(v)}{\log^2 \Delta} \right\} \cap \bigcap \overline{B_w}} = 1 -
  o(e^{-\log \Delta \log \log \Delta}).
\end{equation}

Crucially, the high probability event in \eqref{ready-for-x3} is entirely determined by
Phase I, so all of the choices in Phase II are still independent of it.  Now condition on
Phase I (i.e., $X_2$ and $I_2(k)$ are fixed), and also on the 
event in \eqref{ready-for-x3}.  Perform Phase II. We will show that with high probability the random 
variable $R = X_2 -X_3$ is small.  Observe that since $I_2 \geq I_3$, and we conditioned on Phase I,
\begin{displaymath}
  \E{R} = \sum_{k=1}^r d_t(v, k) \E{I_2(k) - I_3(k)} 
  = \sum_{1 \leq k \leq r, I_2(k) = 1} d_t(v, k) \pr{I_2(k) - I_3(k) = 1}.
\end{displaymath}
Now given that $I_2(k) = 1$, the difference $I_2(k) - I_3(k)$ will be 1 precisely when the
vertex $u \in V_k \cap S$ has one of its (at most $D_t$) neighbors $w$ selected in Phase
II.  For each such neighbor $w$, the probability of its selection in Phase II is $\leq 1/(S_t
\log \Delta)$, so a simple union bound gives 
$\pr{I_2(k) - I_3(k) = 1}\leq \frac{D_t}{S_t \log \Delta}$. Therefore
$$  \E{R} \leq \sum_{1 \leq k \leq r, I_2(k) = 1} d_t(v, k) \frac{D_t}{S_t \log \Delta} 
  = X_2 \frac{D_t}{S_t \log \Delta}  \leq  \frac{X_2}{\log \Delta}  \leq 
  \Theta\left( \frac{d_t(v)}{\log^2 \Delta} \right),$$
since we conditioned on a range for $X_2$.  Next we show that $R$ is 
concentrated.  We
conditioned on $\bigcap \overline{B_w}$, so changing any choice in Phase II can affect $R$
by at most $C \log \Delta$.  Therefore, $R$ is Lipschitz with constant $C \log \Delta$.  It is
also clear that $R$ is 1-certifiable. Since $d_t(v) = \Theta(\Delta)$ and $R \leq X_2 \leq 
\Delta$, by Talagrand's Inequality (Theorem \ref{thm:talagrand}) we have
\begin{eqnarray*}
  \pr{\big|R - \E{R}\big| > \frac{d_t(v)}{\log^2 \Delta}} &\leq& 4 \exp\left\{ 
    - \Theta\left( \left(\frac{d_t(v)}{\log^2 \Delta}\right)^2 \frac{1}{8(C \log \Delta)^2
        \E{R}} \right)
  \right\} \\
  &\leq& \exp\left\{ -\Theta\left( \frac{\Delta}{\log^6 \Delta} \right) \right\} \\
  &\ll& e^{-\log \Delta \log \log \Delta}\,.
\end{eqnarray*}
In particular,
\begin{displaymath}
  \pr{X_2-X_3 > \Theta\left(\frac{d_t(v)}{\log^2 \Delta}\right)} \ll e^{-\log \Delta \log
    \log \Delta}\, .
\end{displaymath}
Therefore, with probability $1-o\left(e^{-\log \Delta \log \log \Delta}\right)$,
$$z_t(v) = X_3 \geq X_2 - \Theta\Big(\frac{d_t(v)}{\log^2 \Delta}\Big)\geq
 \frac{d_t(v)}{\log \Delta} - \Theta\Big(\frac{d_t(v)}{\log^2 \Delta}\Big).$$  
This establishes \eqref{z-goal} and completes the proof.
\hfill $\Box$

\section{Clique-free transversals}
\label{sec:generalize}

In this section, we study sufficient conditions for the existence of a $K_s$-free
transversal in graphs $G$ with maximum degree $\Delta$.  Consider $s$ to be a fixed
parameter, and let $\Delta$ grow.  We will prove that if the local degree is $o(\Delta)$,
then parts of size $(1+o(1)) \frac{\Delta}{s-1}$ are sufficient.

First, let us show that this bound is asymptotically tight via the following construction.
Fix any positive integer $n < \frac{\Delta+1}{s-1}$, and let $G$ be a graph with vertex
set $V = \{1, \ldots, \Delta+1\} \times \{1, \ldots, n\}$.  Let the parts be defined as
$V_i = \{(i, j) : 1 \leq j \leq n \}$, and let $(i, j)$ and $(i', j)$ be adjacent for all
$1 \leq i, i' \leq \Delta+1$.  It is clear that $G$ has maximum degree $\Delta$ and local
degree 1.  We show by contradiction that this graph has no $K_s$-free transversal.
Indeed, if there is such a transversal $T$, then for each $j$, the set of vertices $(i, j)
\in T$ forms a clique and hence has cardinality at most $s-1$.  Yet there are only $n$
possibilities for $j$, so $|T| \leq n(s-1) < \Delta+1$.  This is a contradiction, since
$T$ must have one vertex in each of the $\Delta+1$ parts.  Therefore, parts of size
$\frac{\Delta}{s-1}$ do not guarantee a $K_s$-free transversal.

\vspace{0.2cm}
\noindent {\bf Proof of Theorem \ref{generalize}.}\, Fix $\ep > 0$ and $s \geq 3$.  Let
$G=(V, E)$ be a graph with maximum degree at most $\Delta$ whose vertex set is partitioned
into $r$ parts $V=V_1 \cup \ldots \cup V_r$ of size $|V_i| \geq
(1+\ep)\frac{\Delta}{s-1}$.  Color the vertices of $G$ with $s-1$ colors such that the
number of monochromatic edges is minimal.  Note that for every vertex $v$, there must be a
color $c$ such that the number of neighbors of $v$ which are colored $c$ is at most
$\big\lfloor \frac{\Delta}{s-1} \big\rfloor$.  Hence the minimality of the coloring
implies that $v$ has at most that many neighbors in its own color, or else one could
obtain a better coloring by changing the color of $v$ to $c$.  Now delete all edges whose
endpoints have different colors, and call the resulting graph $G'$.  By the above
argument, the maximum degree in $G'$ is at most $\big\lfloor \frac{\Delta}{s-1}
\big\rfloor$, so $G'$ has an independent transversal $T$ by Theorem \ref{main}.  However,
$T$ is an $(s-1)$-colorable transversal in $G$, and so must be $K_s$-free.  \hfill $\Box$

\vspace{0.2cm} Observe that we did not need the local degree condition until we invoked
Theorem \ref{main}.  If we do not have a local degree condition, we can apply Haxell's
result \cite{Ha-note} instead, which says that parts of size $2\Delta$ guarantee an
independent transversal in graphs with maximum degree $\Delta$.  This immediately implies:

\begin{proposition}
\label{noga}
  Let $G$ be a graph with maximum degree at most $\Delta$ whose vertex set is partitioned
  into $r$ parts $V(G) = V_1 \cup \ldots \cup V_r$ of size $|V_i| \geq
  2 \big\lfloor \frac{\Delta}{s-1} \big\rfloor$.  Then $G$ has a $K_s$-free transversal.
\end{proposition}

\noindent Phrased in terms of the function $p(\Delta, K_s)$ defined in the introduction,
we have
\begin{displaymath}
  p(\Delta, K_s) \leq 2 \left\lfloor \frac{\Delta}{s-1} \right\rfloor,
\end{displaymath}
which is at most twice Szab\'o and Tardos's lower bound (which they conjectured to be tight)
\begin{displaymath}
  p(\Delta, K_s) \geq \frac{s}{s-1} \frac{\Delta}{s-1}.
\end{displaymath}
Note that for $s=3$, it matches their best upper bound, $p(\Delta, K_3) \leq \Delta$,
which they obtain as a consequence of a result on acyclic transversals, i.e., transversals
which have no cycles.  So, this simple approach provides an alternate proof of that upper
bound.  For $s>3$, as far as we know, this proposition gives the current best upper bound.

\section{Concluding Remarks}

\begin{itemize}
\item We proved that if $G$ is a multipartite graph with maximum degree $\Delta$ and local
degree $o(\Delta)$, then parts of size $(1+o(1))\Delta$ will guarantee an independent
transversal.  It is interesting to decide if it is possible to achieve the same result under 
the weaker condition that the number of edges between any pair of distinct parts is 
$o(\Delta^2)$.

\item Let $M=M(\Delta)$ be the smallest integer such that if $G$ is a multipartite graph
  with maximum degree $\Delta$, local degree 1, and parts of size $\Delta+M$, then it has
  an independent transversal. We showed that $M=o(\Delta)$ (in fact, this can be sharpened
  to $\Delta^{1-\ep}$ using our method) and it remains an interesting problem to better
  estimate the function $M(\Delta)$. In particular, an intriguing open question is to
  determine if $M(\Delta)$ is bounded by absolute constant.  Note that a list coloring
  construction of Bohman and Holzman from \cite{BH} implies that $M$ would have to be at
  least $2$, because as mentioned in the introduction, an instance of the list coloring problem
  corresponds to an independent transversal problem with local degree 1.

\item Let $G$ be a graph with maximum degree $\Delta$ whose vertex set is partitioned into
  $r$ equal parts $V(G)=V_1 \cup \ldots \cup V_r$ of size $n$. How large should $n$ be to
  ensure that we can partition the entire graph into a disjoint union of $n$ independent
  transversals? This question is related to the notion of \emph{strong chromatic number},
  see, e.g., \cite{ABZ, Al-strong-chromatic, Ha-strong-chromatic}.  Alon
  \cite{Al-strong-chromatic} proved that for a (large) constant $c$, parts of size
  $n=c\Delta$ are enough.  Haxell \cite{Ha-strong-chromatic} reduced the constant to 3,
  and recently even to $3-\ep$, where $\ep$ can be as large as $1/4$
  \cite{Ha-private-comm}.  It would be very interesting to determine the correct value of
  $c$, which should be at least $2$ because of the construction of Szab\'o and Tardos
  mentioned in the introduction.

  However, if we impose a local degree restriction on $G$, our result suggests that one
  does not need parts of size $2\Delta$. We believe that if $G$ has maximum degree
  $\Delta$ and local degree $o(\Delta)$ then parts as small as $n=(1+o(1))\Delta$ will
  guarantee the existence of $n$ disjoint independent transversals. So far we can only
  prove the much weaker statement that parts of size at least $(2+o(1))\Delta$ are
  sufficient.  This claim follows immediately from our main result together with an
  argument of Aharoni, Berger, and Ziv. In \cite{ABZ} (see Theorem 5.3) they implicitly
  proved that if parts of size at least $f(\Delta)$ imply that every vertex $v$ of $G$ is
  contained in some independent transversal, then parts of size at least $\Delta
  +f(\Delta)$ guarantee the existence of a partition of $G$ into independent transversals.
  Our result certainly implies the former statement with $f(\Delta) = (1+\ep)\Delta$.
  Indeed, for any given vertex $v$, the local degree is $o(\Delta)$, so we can delete
  $o(\Delta)$ neighbors of $v$ from every part.  Then $v$ becomes isolated from rest of
  the graph.  However, the part sizes are still at least $(1+\ep-o(1))\Delta$ so by
  Theorem \ref{main} we can find an independent transversal among the parts not containing
  $v$, and then add $v$.
\end{itemize}

\vspace{0.2cm}
\noindent {\bf Acknowledgment.}\, The authors would like to thank Noga Alon, whose
suggestion simplified their original proofs of Theorem \ref{generalize} and Proposition
\ref{noga}, and the referees for careful reading of this manuscript.

\end{document}